\input amstex
\documentstyle{amsppt}
\topmatter
\title
Construction of Vector Valued Modular Forms from Jacobi Forms
\endtitle
\rightheadtext{Construction of Modular Forms from Jacobi Forms}
\author  Jae-Hyun Yang
\endauthor
\magnification=\magstep 1 \pagewidth{12.5cm} \pageheight{18.20cm}
\baselineskip =7mm

%Use \endgraf to indicate new paragraph
%\thanks will become a 1st page footnote
\thanks{This work was in part supported by TGRC-KOSEF and
Max-Planck-Institut f{\" u}r Mathematik.\endgraf
1991 Mathematics Subject Classification. Primary 11F30, 11F55. }
\endthanks
\abstract{We give a geometrical construction of the canonical automorphic
factor for the Jacobi group and construct new vector valued modular forms
from Jacobi forms by differentiating them with respect to toroidal
variables and then evaluating at zero.}
\endabstract
\endtopmatter
\document
\NoBlackBoxes

\def\G{\Gamma}

\def\l{\lambda}

\def\lb{\lbrace}

\def\rb{\rbrace}

\def\s{\sigma}

\def\lrt{\longrightarrow}
\def\lmt{\longmapsto}

\def\M{{\cal M}}
\def\Box{\square}

\def\s{\sigma}

\def\lrt{\longrightarrow}
\def\lmt{\longmapsto}

\def\M{{\Cal M}}
\def\Hom{\text {Hom}}

\def\cH{{\Cal H}}
\def\cM{{\Cal M}}
\def\BZ{\Bbb Z}

\def\BC{\Bbb C}
\def\BR{\Bbb R}
\def\G{\Gamma}
\head 1 \ Introduction %bold
\endhead               %don't type final punctuation
For given two fixed positive integers $n$ and $m$, we let
$$H_n:=\{\,Z\in {\BC}^{(n,n)}\,\vert\ Z=\,^tZ,\ \ \text {Im}\,Z\,>\,0\ \}$$
be the Siegel upper half plane of degree $n$ and let $\Gamma_n$ be the
Siegel modular group of degree $n$. Let
$${\Cal P}_{m,n}:={\BC}[W_{11},\cdots,W_{mn}],\ \ \ W=(W_{kl})\in
{\BC}^{(m,n)}$$
be the ring of polynomial functions on ${\BC}^{(m,n)}.$ Here ${\BC}^{(n,n)}\,
(resp.\,{\BC}^{(m,n)})$ denotes the space of all complex $n\times n\,(resp.\,
m\times n)$-matrices (\,see notations below\,). For any homogeneous
polynomial $P\in {\Cal P}_{m,n},$ we define the differential
operator $P(\partial_W)$ on ${\BC}^{(m,n)}$ as follows:
$$P(\partial_W):=P\left({{\partial}\over {\partial W_{11}}},\cdots,
{{\partial}\over {\partial W_{mn}}}\right).$$
In this paper, the author proves that if $P$ is a
{\it homogeneous\ pluriharmonic}
polynomial in ${\Cal P}_{m,n}$ and $f\in J_{\rho,\M}(\G_n)$
(see Definition 3.1) is a {\it
Jacobi\ form} of index $\M$ with respect to a rational representation
$\rho$ of the general group $GL(n,\BC),$ then the following function
$$P(\partial_W)\,f(Z,W)\Big\vert_{W=0}$$
yields a vector valued modular form with respect to a new rational
representation of $GL(n,\BC).$ For a precise detail, we refer to Definition
5.1 and Main Theorem in section 5. In [M-N-N] (cf. pp 147-156), the authors
proved the similiar result for theta functions. Our result is a
generalization of their result because theta functions are special examples
of Jacobi forms.
\par
\indent
This paper is organized as follows.
In section 2, we provide a geometrical construction of the canonical
automorphic factor for the Jacobi group.
In section 3, we review Jacobi forms and establish the notations.
In section 4, we review pluriharmonic polynomials and obtain some
properties to be used in the subsequent sections. In section 5,
we shall prove the main theorem.
In section 6, we obtain two identities by applying the main theorem to
Jacobi forms.
\par\smallpagebreak\noindent
{\smc
Notations:} \ \ We denote by $\BZ,\,\BR$ and $\BC$ the ring of integers,
the field of
real numbers, and the field of complex numbers respectively.
$\G_n:=Sp(n,\BZ)$ denotes the Siegel modular
group of degree $n$.
The symbol
``:='' means that the expression on the right is the definition of that on the
left. We denotes by ${\BZ}^+$ the set of
all positive integers. $F^{(k,l)}$ denotes
the set of all $k\times l$ matrices with entries in a commutative ring $F$.
For a square matrix $A\in F^{(k,k)}$ of degree $k$, $\sigma(A)$ denotes the
trace of $A$.
For $A\in F^{(k,l)}$ and $B\in F^{(k,k)},$ we set $B[A]=\,^t\!ABA.$ For
any $M\in F^{(k,l)},\ ^t\!M$ denotes the transpose matrix of $M$. $E_n$ denotes
the identity matrix of degree $n$.\par\medpagebreak
\vskip 0.3cm
\head 2 \ The\ Canonical\ Automorphic\ Factor\ for\ the\ Jacobi
\ Group            %bold
\endhead           %don't type final punctuation
Let $m$ and $n$ be two fixed positive integers. It is well known that the
automorphism group $Aut\,(H_{m+n})$ of the Siegel upper half plane of degree
$m+n$ is given by
$$Aut\,(H_{m+n})=Sp(m+n,\BR)/\{\pm E_{m+n}\}.$$
We observe that $H_n$ is a rational boundary of $H_{m+n}$\,(cf.\,[N]).
The normalizer $N(H_n):=\{ {\tilde \sigma}\in Aut(H_{m+n})\,:\
{\tilde \sigma}(H_n)\subset H_n\ \}$ of $H_n$ is given by
$$N(H_n)=P(H_n)/\{ \pm E_{m+n}\},$$
where
$$\align
P(H_n):&=\left\{ g\in Sp(m+n,\BR)\ :\ g(H_n)\subset H_n\ \right\} \\
&=\left\{ [\sigma,u,(\l,\mu,\kappa)]\in Sp(m+n,\BR)\right\}.
\endalign$$
Here we put
$$[\sigma,u,(\l,\mu,\kappa)]:=
\left( \matrix A & 0 & B & A\,{^t\mu}-B\,{^t\!\l}\\
u\l & u & u\mu & u\kappa\\
C & 0 & D & C\,{^t\!\mu}-D\,{^t\!\l}\\
0 & 0 & 0 & ^t\!u^{-1}\endmatrix\right) ,$$
where $\sigma=
\left(\matrix A & B \\ C & D\endmatrix\right)
\in Sp(n,\BR),\ u\in GL(m,\BR),
\ \l,\mu\in \BR^{(m,n)}$ and $\kappa\in \BR^{(m,m)}.$\par
If
$\left(\matrix Z & ^t\!W\\ W & T\endmatrix\right)
\in H_{m+n}$ with $Z\in H_n,\,
W\in \BR^{(m,n)}$ and $T\in H_m,$ we simply write
$$(Z,W,T):=
\left(\matrix Z & ^tW\\ W & T\endmatrix\right).$$
We denote the symplectic action of $N(H_n)$ on $(Z,W,T)$ by
$$g\cdot (Z,W,T):=({\tilde Z},{\tilde W},{\tilde T}),\ \ \ \
g\in N(H_n).$$
It is easy to see that $({\tilde Z},{\tilde W},{\tilde T})$ is of the form
$$\align
{\tilde Z}&=\sigma_g(Z),\\
{\tilde W}&=a(g;Z)(W)+b(g;Z),\\
{\tilde T}&=m_g(T)+c(g;Z,W),
\endalign$$
where $\sigma_g\in Aut\,(H_n),\ m_g\in Aut\,({\Cal P}_m),$
$$\align
a(g;\cdot)\,&:\,H_n\lrt GL(\BC^{(m,n)})\ \ \ holomorphic,\\
b(g;\cdot)\,&:\,H_n\lrt \BC^{(m,n)}\ \ \ \ holomorphic,\\
c(g;\cdot,\cdot)\,&:\,H_n\times \BC^{(m,n)}\lrt H_m\ \ \ holomorphic.
\endalign$$
Here ${\Cal P}_m:=\{\,Y\in \BR^{(m,m)}\,|\ Y=\,^tY>0\ \}$ is an open convex
cone in $\BR^{{m(m+1)}\over 2}$ and we set
$$Aut\,({\Cal P}_m):=
\left\{ \,\xi\in GL(\BC^{(m,m)})\,\bigg| \ \xi({\Cal P}_m)
={\Cal P}_m\ \right\}.$$
\noindent
{\smc Remark\ 2.1.}\ In [PS], Piateski-Sharpiro called the mapping
$(Z,W,T)\longmapsto ({\tilde Z},{\tilde W},{\tilde T})$ a {\it quasilinear}
transformation.\par\medpagebreak
From now on, we set
$$H_{n,m}:=H_n\times \BC^{(m,n)}.$$
We observe that $g=[\sigma,u,(\l,\mu,\kappa)](mod\,\{\pm E_{m+n}\})\in
N(H_n)$ acts on $H_{n,m}$ by
$$(Z,W)\lmt (\sigma_g(Z),a(g;Z)(W)+b(g;Z)).$$
The subgroup of $N(H_n)$ consisting of elements $g=[\sigma,u,(\l,\mu,\kappa)]
(mod\,\{\pm E_{m+n}\})$ with the property
$$m_g=Identity \ \ \ \ \ on\ H_m$$
is called the {\it Jacobi\ group}, denoted by $G^J$. It follows immediately
from the definition that
$$G^J=\{[\sigma,E_m,(\l,\mu,\kappa)]\in P(H_n)\}.$$
It is easy to see that $G^J$ is the semidirect product of $Sp(n,\BR)$ and
$H_{\BR}^{(n,m)}$, where
$$H_{\BR}^{(n,m)}:=\left\{ [E_n,E_m,(\l,\mu,\kappa)]\in P(H_n)\right\}$$
is the nilpotent 2-step subgroup of $P(H_n)$, called the {\it Heisenberg\
group}. For some results on $H_{\BR}^{(n,m)},$ we refer to [Y1]-[Y2].
\par\medpagebreak\indent
Now we consider another subgroup ${\tilde G}$ of $G^J.$ By
the definition, ${\tilde G}$ consists of elements of $G^J$ whose action is of
the following form:
$$(Z,W,T)\lmt (\s_g(Z),a(g;Z)(W),T+c(g;Z,W)),\ \ \ c(g;Z,0)=0.$$
\noindent
{\smc Lemma\ 2.2.}\ The map
$$J\,:\,{\tilde G}\times H_n\lrt GL(\BC^{(m,n)})$$
defined by
$$J({\tilde \sigma},Z):=a({\tilde \sigma};Z),\ \ {\tilde \sigma}\in
{\tilde G},\ \ Z\in H_n$$
is a factor of automorphy for ${\tilde G}.$
\vskip 0.15cm
\noindent
{\it Proof.} It is easy to prove it. We leave its proof to the reader.
\hfill $\Box$\par\medpagebreak\noindent
We note that the mapping
$$A(g,(Z,W)):=c(g;Z,W),\ \ g\in G^J,\ (Z,W)\in H_{n,m}
\tag 2.1$$
is a summand of automorphy, i.e.,
$$A(g_1g_2,(Z,W))=A(g_1,g_2\cdot (Z,W))+A(g_2,(Z,W)),\tag 2.2$$
where $g_1,\, g_2\in G^J$ and $(Z,W)\in H_{n,m}.$
We let
$$K_{\BC}\subset GL(\BC^{(m,n)})$$
be the complex Lie group generated by the linear mapping
$$\left\{\ a(g;Z)\,:\ g\in G^J\ \right\}.$$
Then $K_{\BC}$ is isomorphic to $GL(n,\BC).$\par\medpagebreak
\noindent
{\smc Lemma\ 2.3.}\ Let
$$\rho : GL(n,\BC)\lrt GL(V_{\rho})$$
be a finite dimensional holomorphic representation of $GL(n,\BC)$ on a
finite dimensional complex vector space $V_{\rho}$ and let
$\chi\, :\, \BC^{(m,n)}\lrt \BC^{\times}$
be a character on the additive group $\BC^{(m,m)}$. Then the mapping
\,$J_{\rho}\,:\,{\tilde G}\times H_n\lrt GL(V_{\rho})$\,
defined by
$$J_{\rho}({\tilde \sigma},Z):=\rho(J({\tilde \sigma},Z)),\ \ \
{\tilde \sigma}\in {\tilde G},\ \ Z\in H_n$$
is a factor of automorphy for ${\tilde G}$. Furthermore the mapping
$$J_{\chi,\rho}(g,(Z,W)):=\chi(c(g;Z,W))\,\rho(a(g;Z)),\ \ g\in G^J$$
is a factor of automorphy for the Jacobi group $G^J$ with respect to
$\chi$ and $\rho.$
\vskip 0.1cm
\noindent
{\it Proof.}\
The proof of this first statement is obvious. The proof of the second
statement follows immediately from the fact that $A(g,(Z,W)):=c(g;Z,W)$ is
a summand of automorphy (cf. (2.1) and (2.2)) and that $J_{\rho}$ is a
factor of automorphy for ${\tilde G}.$ \hfill $\Box$
\par\medpagebreak
\noindent
{\smc Definition\ 2.4.}\ $J_{\rho}$ and $J_{\chi,\rho}$ are called the
{\it canonical\ automorphic\ factor} for ${\tilde G}$ with respect to
$\rho$ and the {\it canonical\ automorphic\ factor} for $G^J$ with respect
to $\chi$ and $\rho$ respectively.
\vskip 0.15cm
\head 3 \ Jacobi\ Forms %bold
\endhead                %don't type final punctuation
\vskip 0.2cm
In this section, we establish the notations and define the concept of
Jacobi forms.\par\smallpagebreak
Let
$$Sp(n,\BR)=\lb M\in \BR^{(2n,2n)}\ \vert \ ^t\!MJ_nM= J_n\ \rb$$
be the symplectic group of degree $n$, where
$$J_n:=\left(\matrix 0&E_n\\
                   -E_n&0\endmatrix\right).$$
It is easy to see that $Sp(n,\BR)$ acts on $H_n$ transitively
by
$$M<Z>:=(AZ+B)(CZ+D)^{-1},$$
where $M=\left(\matrix A&B\\ C&D\endmatrix\right)\in Sp(n,\BR)$ and $Z\in H_n.$
\par\smallpagebreak
For two positive integers $n$ and $m$, we recall that
the Jacobi group $G^J:=Sp(n,\BR)\ltimes H_{\BR}^{(n,m)}$ is
the semidirect product of the symplectic group $Sp(n,\BR)$ and
the Heisenberg group $H_{\BR}^{(n,m)}$
endowed with the following multiplication law
$$
(M,(\lambda,\mu,\kappa))\cdot(M',(\lambda',\mu',\kappa'))
:=\, (MM',(\tilde{\lambda}+\lambda',\tilde{\mu}+
\mu', \kappa+\kappa'+\tilde{\lambda}\,^t\!\mu'
-\tilde{\mu}\,^t\!\lambda'))$$
with $M,M'\in Sp(n,\BR), (\lambda,\mu,\kappa),\,(\lambda',\mu',\kappa')
\in H_{\BR}^{(n,m)}$
and $(\tilde{\lambda},\tilde{\mu}):=(\lambda,\mu)M'$.
It is easy to see that $G^J$
acts on $H_{n,m}:=H_n\times \BC^{(m,n)}$ transitively by
$$(M,(\lambda,\mu,\kappa))\cdot (Z,W):=(M<Z>,(W+\lambda Z+\mu)
(CZ+D)^{-1}),
\tag 3.1$$
where $M=\left(\matrix A&B\\ C&D\endmatrix\right)
\in Sp(n,\BR),\ (\lambda,\mu,
\kappa)\in H_{\BR}^{(n,m)}$ and $(Z,W)\in H_{n,m}.$
\vskip 0.3cm
\indent
Let $\rho$ be a rational representation of $GL(n,\BC)$ on a finite dimensional
complex vector space $V_{\rho}.$ Let ${\Cal M}\in \BR^{(m,m)}$ be a symmetric
half-integral semi-positive definite matrix of degree $m$.
Let $C^{\infty}(H_{n,m},V_{\rho})$ be the algebra of all
$C^{\infty}$ functions on $H_{n,m}$
with values in $V_{\rho}.$ For $f\in C^{\infty}(H_{n,m},
V_{\rho}),$
we define
$$\align
  & (f|_{\rho,{\Cal M}}[(M,(\lambda,\mu,\kappa))])(Z,W)\\
(3.2)\quad\quad:=
\,& e^{-2\pi i\sigma({\Cal M}[W+\lambda Z+\mu](CZ+D)^{-1}C)}
\times
e^{2\pi i\sigma({\Cal M}(\lambda Z^t\!\lambda+2\lambda^t\!W+(\kappa+
\mu^t\!\lambda)))} \\
&\times\rho(CZ+D)^{-1}f(M<Z>,(W+\lambda Z+\mu)(CZ+D)^{-1}),
\hskip 3.8cm
\endalign$$
where $M=\left(\matrix A&B\\ C&D\endmatrix\right)\in Sp(n,\BR),\
(\l,\mu,\kappa)\in H_{\BR}^{(n,m)}$
and $(Z,W)\in H_{n,m}.$\par\medpagebreak
\noindent
{\smc Definition\ 3.1.}\ \ Let $\rho$ and $\M$ be as above. Let
$$H_{\BZ}^{(n,m)}:=\lb \, (\l,\mu,\kappa)\in H_{\BR}^{(n,m)}\, \vert
\, \l,\mu\in \BZ^{(m,n)},\ \kappa\in \BZ^{(m,m)}\, \rb.$$
A {\it Jacobi\ form} of index $\M$ with respect to $\rho$
on $\G_n$ is a holomorphic
function $f\in C^{\infty}(H_{n,m},V_{\rho})$ satisfying the
following conditions (A) and (B):\par\medpagebreak
\noindent
(A)\ \ \ \
$f|_{\rho,{\Cal M}}[\tilde{\gamma}] = f$ for all $\tilde{\gamma}\in
\Gamma^J_n := \G_n \ltimes H_{\BZ}^{(n,m)}$.\par\medpagebreak
\noindent
(B)\ \ \ $f$ has a Fourier expansion of the following form :
$$f(Z,W) = \sum\limits_{T\ge0\atop \text {half-integral}}
\sum\limits_{R\in \BZ^{(n,m)}}
c(T,R)\cdot e^{{2\pi i}\,\sigma(TZ)}\cdot
e^{2\pi i\sigma(RW)}$$
with $c(T,R)\ne 0$ only if $\left(\matrix
T & \frac 12R\\ \frac 12\,^t\!R&{\Cal M}\endmatrix\right)
\ge 0$.
\par
\indent
If $n\geq 2,$ the condition (B) is superfluous by
K{\" o}cher principle\,(\,cf.\,
[Z] Lemma 1.6). We denote by $J_{\rho,\M}(\G_n)$
the vector space of all
Jacobi forms of index $\M$
with respect to $\rho$ on $\G_n$.
Ziegler\,(\,cf.\,[Z] Theorem 1.8 or [E-Z] Theorem 1.1\,)
proves that the vector space $J_{\rho,\M}(\G_n)$ is finite dimensional.
For more results on Jacobi forms with $n>1$ and $m>1$, we refer to
[Y3]-[Y6] and [Z].
\vskip 0.2cm
\head 4 \ Pluriharmonic\ Polynomials   %bold
\endhead                               %don't type final punctuation
\vskip 0.2cm
We review pluriharmonic polynomials of matrix arguments and collect some
properties to be used in the next section
\,(\,cf. [K-V] and [M-N-N]\,).\par\smallpagebreak
\indent
Let $n$ and $m$ be two positive integers and let ${\Cal P}_{m,n}:=
\BC[W_{11},W_{12},\cdots,W_{mn}]$ be the ring of complex
valued polynomials on $\BC^{(m,n)}.$ For any homogeneous polynomial
$P\in {\Cal P}_{m,n},$ we put
$$P(\partial_W):=P\left({{\partial}\over {\partial W_{11}}},\cdots,
{{\partial}\over {\partial W_{mn}}}\right). \tag 4.1$$
Let $S$ be a positive definite symmetric rational matrix of degree $m$.
Let $T:=(t_{pq})$ be the inverse of $S$. For each $i,j$ with $1\leq i,j
\leq n,$ we denote by $\triangle_{i,j}$ the following differential operator
$$\triangle_{i,j}:=\sum_{p,q=1}^m\,t_{pq}
{{\partial^2}\over {\partial W_{pi}\,\partial W_{qj}}},\ \ \ 1\leq i,j\leq
n. \tag 4.2$$
A polynomial $P$ on $\BC^{(m,n)}$ is said to be {\it harmonic} with respect
to $S$ if
$$\sum_{i=1}^n\, \triangle_{i,i}P=0.\tag 4.3$$
A polynomial $P$ on $\BC^{(m,n)}$ is called {\it pluriharmonic} with respect
to $S$ if
$$\triangle_{i,j}P=0,\ \ \ \ 1\leq i,j\leq n.\tag 4.4$$
If there is no confusion, we just write harmonic or pluriharmonic instead of
harmonic or pluriharmonic with respect to $S$. Obviously a pluriharmonic
polynomial is harmonic. We denote by ${\Cal H}_{m,n}$ the space of all
pluriharmonic polynomials on ${\BC}^{(m,n)}.$ The ring ${\Cal P}_{m,n}$ of
polynomials on ${\BC}^{(m,n)}$ has a symmetric nondegenerate bilinear form
$<P,\,Q>:=(P(\partial_W)Q)(0)$ for $P,\,Q\in {\Cal P}_{m,n}.$ It is easy to
check that $<\,,\,>$ satisfies
$$<P,\,QR>=<Q(\partial_W)P,\,R>,\ \ \ \ P,Q,R \in {\Cal P}_{m,n}.\tag 4.5$$
\vskip 0.2cm
\noindent
{\smc Lemma 4.1.}\ \ ${\Cal H}_{m,n}$ is invariant under the action of
$GL(n,\BC)\times O(S)$ given by
$$((A,B),\,P(W))\longmapsto P(\,^tBWA),\ \ A\in GL(n,\BC),\ B\in O(S).
\tag 4.6$$
Here $O(S):=\{\,B\in GL(m,\BC)\,\vert\,\ ^tBSB=S\ \}$ denotes
the orthogonal group
of the quadratic form $S$.\par\medpagebreak
\noindent
{\it Proof.}\ See corollary 9.11 in [M-N-N].
\hfill $\Box$ \par\smallpagebreak
\noindent
{\smc Remark\ 4.2.}\ In [K-V], Kashiwara and Vergne investigated an
irreducible decomposition of the space of complex pluriharmonic polynomials
defined on ${\BC}^{(m,n)}$ under the action of (4.6).
They showed that each
irreducible component $\tau\otimes \lambda$
occurring in the decomposition of ${\Cal H}_{m,n}$
under the action (4.6) has multiplicity one and the irreducible
representation $\tau$ of $GL(n,\BC)$ is determined uniquely by the irreducible
representation of $O(S).$
\vskip 0.2cm
\noindent
{\smc Lemma\ 4.3.}\ \ If $P$ is pluriharmonic, then we have
$$P(\partial_W)\,e^{\sigma(WC\,^tWS^{-1})}=P(2S^{-1}WC)\,
e^{\sigma(WC\,^tWS^{-1})}$$
for all complex {\it symmetric} matrix $C\in {\BC}^{(n,n)}$ of degree $n.$
We recall that $\sigma(A)$ denotes the trace of a square matrix.
\vskip 0.2cm
\noindent
{\it Proof.}\ \ We set $h(W):=\sigma(WC\,^tWS^{-1}).$ We observe that
$h(\partial_W)P=0.$ Indeed,
$$\align
h(W)&=\sum_{i,k,l,m}W_{ik}c_{kl}W_{ml}t_{mi}\\
    &=\sum_{k,l}c_{kl}\left(\sum_{i,m}t_{mi}W_{ml}W_{ik}\right) \\
    &=\sum_{k,l}c_{kl}h_{lk}.
\endalign$$
Thus $h(\partial_W)P=\sum_{k,l}c_{kl}\left(h_{lk}(\partial_W)P\right)=
\sum_{k,l}\,c_{kl}\,\triangle_{l,k}P=0.$ We put $\varphi(W):=
e^{h(W)}.$ Then $\,f(W):=\varphi(W+A)=\varphi(W)\varphi(A)\eta(W),$
where $A\in {\BC}^{(m,n)}$ and $\eta(W):=e^{\sigma(2WC{^t\!A}S^{-1})}.$
$$\align
P(\partial_W)\varphi(W)\Big\vert_{W=A}
&=P(\partial_W)f(W)
\Big\vert_{W=0}\\
&=\varphi(A)\left(P(\partial_W)\varphi(W)\eta(W)\right)\Big\vert_{W=0}\\
&=\varphi(A)P(\partial_W)\eta(W)\Big\vert_{W=0}.
\endalign$$
Indeed, since $h(\partial_W)P=0,$ we have
$$\align
P(\partial_W)\left(\varphi(W)\eta(W)\right)\Big\vert_{W=0}
&=<P,\,\varphi\cdot\eta>=<\varphi(\partial_W)P,\,\eta>\\
&=\sum_{n=0}^{\infty}{1\over {n!}}<h^n(\partial_W)P,\,\eta>\\
&=<P,\eta>=P(\partial_W)\eta(W)\Big\vert_{W=0}.
\endalign$$
By an easy computation, we obtain
$$P(\partial_W)\,\eta(W)=P(2S^{-1}AC)\,\eta(W).$$
Finally, we have
$$P(\partial_W)\varphi(W)\Big\vert_{W=A}=
\varphi(A)\cdot P(2S^{-1}AC)\,\eta(0).$$
Hence we obtain the desired result.\hfill $\Box$
\vskip 0.3cm
\head 5\ Proof\ of\ Main\ Theorem   %bold
\endhead                            %don't type final punctuation
\vskip 2mm
Throughout this section we fix a rational representation
$\rho$ of $GL(n,\BC)$
on a finite dimensional complex vector space $V_\rho$ and
a positive definite symmetric, half-integral
matrix ${\Cal M}$ of degree $m$ once and for all.\par
\indent
We set $S:=(2\cM)^{-1}$. As in the previous section, we
denote by $\cH_{m,n}$ the vector space of all pluriharmonic polynomials with
respect to $S$ on $\BC^{(m,n)}$. According to Lemma 4.1, there exists an
irreducible subspace $V_\tau(\ne 0)$ invariant under the
action of $GL(n,\BC)$
given by (4.6). We denote this represetation by $\tau$. Then we have
$$(\tau(A)P)(W)=P(WA),\quad A\in GL(n,\BC),\ \ P\in V_\tau,\ \
W\in \BC^{(m,n)}.
\tag 5.1$$
The action $\hat{\tau}$ of $GL(n,\BC)$ on $V^*_\tau$ is defined by
$$(\hat{\tau}(A)^{-1}\zeta)(P) := \zeta(\tau(\,^t\!A^{-1})P),\tag 5.2$$
where $A\in GL(n,\BC),\ \zeta\in V^*_\tau$ and $P\in V_\tau$.
\par\medpagebreak
\noindent
{\smc Definition 5.1.}\ \ Let $f\in J_{\rho,\cM}(\Gamma_n)$ be a Jacobi form
of index ${\Cal M}$ with respect to $\rho$  on $\Gamma_n$.
Let $P\in V_\tau$ be a
homogeneous pluriharmonic polynomial. We put
$$f_P(Z) := P(\partial_W)f(Z,W)\Big\vert_{W=0},
\quad Z\in H_n,\ \ W\in \BC^{(m,n)}.
\tag 5.3$$
Now we define the mapping
$$f_\tau : H_n\longrightarrow V_\tau^*\otimes V_\rho$$
by
$$(f_\tau(Z))(P) := f_P(Z),\quad Z\in H_n,\ \ P\in V_\tau.\tag 5.4$$
\noindent
{\smc Definition 5.2.}\ \ A holomorphic function $f:H_n\rightarrow V_\rho$ is
called a {\it modular form of type $\rho$ on $\Gamma_n$} if
$$f(M<Z>)=\rho(CZ+D)f(Z),\quad Z\in H_n$$
for all
$M=\left(\matrix A&B\\ C&D
\endmatrix\right) \in\Gamma_n$. If $n=1$, the additional cuspidal
condition will be added. We denote by $[\Gamma_n,\rho]$ the vector space of all
modular forms of type $\rho$ on $\Gamma_n$.
\vskip 0.3cm
\noindent
{\bf Main Theorem.} \ \ Let $\tau$ and $\hat{\tau}$ be as before.
Let $f\in J_{\rho,
\cM}(\Gamma_n)$ be a Jacobi form. Then $f_\tau(Z)$ is a modular form of type
$\hat{\tau}\otimes\rho$, i.e., $f_\tau\in[\Gamma_n,\hat{\tau}\otimes\rho]$.
\vskip 0.2cm
\noindent
{\it Proof.}\ \ \ Let
$$f(Z,W)=\sum\limits_{T,R}\,c(T,R)\,e^{2\pi i\sigma(TZ)}
\cdot e^{2\pi i\sigma(RW)}$$
be a Fourier expansion of $f(Z,W)$. Then we have
$$P(\partial_W)f(Z,W)=\sum\limits_{T,R}\,P(2\pi i\,^t\!R)\cdot c(T,R)\cdot
e^{2\pi i\sigma(TZ+RW)}$$
and
$$f_P(Z):=P(\partial_W)f(Z,W)\Big\vert_{W=0}=\sum_{T,R}
\, P(2\pi i\, ^t\!R)\cdot
e^{2\pi i\sigma(TZ)}\cdot c(T,R)\tag 5.5$$
Since $f\in J_{\rho,\cM}(\Gamma_n)$, we have the following transformation
law
$$f(\!M\!<\!Z\!>,W(CZ+D)^{-1})=
e^{2\pi i\sigma(\cM W(CZ+D)^{-1}C\,^t\!W)}\cdot
\rho(CZ+D)f(Z,W)  \tag 5.6$$
for all $M=
\left(\matrix A&B\\ C&D\endmatrix\right) \in \Gamma_n$.
Applying $P(\partial_W)$ to (5.6), according to Lemma 4.3, we have
$$\align
&P(\partial_W)f(M<Z>,W(CZ+D)^{-1})\\
=\,\,& P(4\pi i {\Cal M}W(CZ+D)^{-1}C)
\,e^{2\pi i\sigma(\cM W(CZ+D)^{-1}C\,^t\!W)}\\
&\times \rho(CZ+D)f(Z,W)\,+\,h(Z,W)\,+\,
e^{2\pi i\sigma(\cM W(CZ+D)^{-1}C\,^t\!W)}\\
&\times \sum_{T,R}P(2\pi i\,^t\!R)\cdot\rho(CZ+D)c(T,R)\cdot e^{2\pi i\sigma
(TZ+RW)},
\endalign$$
where $h(Z,W)$ is a $V_{\rho}$-valued function
on $H_{n,m}$ whose restriction
to $W=0$ vanishes.
Here we used the fact that $(CZ+D)^{-1}C$ is a complex symmetric matrix of
degree $n$ and Lemma 4.3.
If we evaluate this at $W=0$, $P$ being homogeneous, we have
$$\align
(5.7)\hskip 25mm & P(\partial_W)f(M<Z>,W(CZ+D)^{-1})\Big\vert_{W=0}\\
=\,&\sum\limits_{T,R}P(2\pi i\,^t\!R)
\cdot e^{2\pi i\sigma(TZ)}\cdot \rho(CZ+D) c(T,R).
\hskip 5cm
\endalign$$
On the other hand,
$$\align
&P(\partial_W)f(M<Z>,W(CZ+D)^{-1})\Big\vert_{W=0}\\
=\,&P(\partial_W)\sum\limits_{T,R}\,c(T,R)e^{2\pi i\sigma(T\cdot M<Z>)}\cdot
e^{2\pi i\sigma(RW(CZ+D)^{-1})}\Big\vert_{W=0}\\
=\,&\sum\limits_{T,R}\,P(2\pi i\,^t\!R\,^t\!(CZ+D)^{-1})\cdot e^{2\pi i\sigma(T
\cdot M<Z>)}\cdot c(T,R).
\endalign$$
Thus according to (5.7), we have
$$\align
(5.8)\quad\quad\quad & \sum\limits_{T,R}\tilde{P}(2\pi i\,^t\!R)\cdot e^{2\pi i\sigma
(T\cdot M<Z>)}\cdot c(T,R)\\
=\,&\sum\limits_{T,R}P(2\pi i\,^t\!R)\cdot e^{2\pi i\sigma(TZ)}\cdot
\rho(CZ+D)c(T,R),\hskip 5cm
\endalign$$
where $\tilde{P}(W):=P(W\,^t\!(CZ+D)^{-1})$. By (5.5), (5.8) implies
$$f_{\tilde{P}}(M<Z>) = \rho(CZ+D)f_P(Z),\tag 5.9$$
that is,
$$(f_\tau(M<Z>))(\tilde{P}) = \rho(CZ+D)f_\tau(Z)(P).\tag 5.10$$
Since $\tilde{P} = \tau(\,^t\!(CZ+D)^{-1})P$, we have from (5.9)
$$\align
&\Big((\hat{\tau}^{-1}\otimes 1_{V_\rho})(CZ+D)f_\tau(M<Z>)\Big)(P)\\
=\,&\Big((1_{{V_\tau}^*}\otimes \rho)(CZ+D) f_\tau(Z)\Big)(P),
\endalign$$
where $1_{{V_\tau}^*}\,(resp.\,V_{\rho})$ denotes the trivial representation
of $GL(n,\BC)$ on $V_{\tau}^*\,(\,resp.\,V_\rho).$
Hence we obtain
$$f_\tau(M<Z>) = (\hat{\tau}\otimes\rho)(CZ+D)f_\tau(Z)\tag 5.11$$
for all $M=
\left(\matrix A&B\\ C&D
\endmatrix\right) \in \Gamma_n$. Therefore $f_\tau$ is a
$\text {Hom}(V_\tau,V_\rho)$-valued
modular form of type $\hat{\tau}\otimes\rho$.
\hfill $\Box$
\vskip 0.31cm
\head 6\ Applications   %bold
\endhead                %don't type final punctuation
\vskip 2mm
In this final section, we
obtain important identites by applying the main theorem
to two special Jacobi forms.
\par\noindent
(I) Let $S\in\BZ^{(2k,2k)}$ be a positive definite
symmetirc, unimodular even matrix of
degree $2k$. We choose an integral matrix $c\in\BZ^{(2k,m)}$
such that $^t\!c
Sc$ is positive definite.
We consider the following theta series
$$\theta_{S,c}(Z,W) :=\sum\limits_{\lambda\in\BZ^{(2k,n)}}e^{\pi i\sigma(S(
\lambda Z\,^t\!\lambda+2\lambda\,^t\!(cW)))}.$$
Then $\theta_{S,c}\in J_{k,\cM}(\Gamma_n)$ with $\cM:=\frac 12\,^t\!cSc$\,
(cf.\,[Z],\,p.\,212).
We write $f(Z,W):=\theta_{S,c}(Z,W)$.
Then by Main Theorem, $f_\tau$ is a $\text {Hom}(V_\tau,\BC)$-valued
modular form of
type $\hat{\tau}\otimes\det^k$. Furthermore, according to (5.9), for any
homogeneous pluriharmonic $P$ with respect to $(2\cM)^{-1}=(\,^t\!cSc)^{-1}$,
we obtain the following identity
$$\align
&\sum\limits_{\lambda\in\BZ^{(2k,n)}}
P(2\pi i\,^t\!cS\lambda\,^t\!(CZ+D)^{-1})
\cdot e^{\pi i\sigma(S\lambda(AZ+B)(CZ+D)^{-1}\,^t\!\lambda)}\\
=\,\,&\{\det\,(CZ+D)\}^k
\sum\limits_{\lambda\in \BZ^{(2k,n)}}
P(2\pi i\,^t\!cS\lambda)\cdot e^{\pi i\sigma(
S\lambda Z\,^t\!\lambda)}
\endalign$$
for all $M=
\left(\matrix A&B\\ C&D\endmatrix\right)
\in\Gamma_n$ and $Z\in H_n$.
\vskip 0.2cm\noindent
(II) In [Z], Ziegler defined the Eisenstein series $E_{k,{\Cal M}}^{(n)}
(Z,W)$ of Siegel type. Let ${\Cal M}$ be a half integral positive definite
symmetric matrix of degree $m$ and let $k\in \BZ^+.$ We set
$$\Gamma_{n,0}:=\left\{
\left(\matrix A & B\\ C & D\endmatrix\right)\in \Gamma_n\Bigg|
\ C=0 \ \right\}.$$
Let ${\Cal R}$ be a complete system of representatives of the cosets
$\Gamma_{n,0}\backslash \Gamma_n$ and $\Lambda$ be a complete system of
representatives of the cosets $\BZ^{(m,n)}/(ker\,({\Cal M})\cap \BZ^{(m,n)}),$
where $ker\,({\Cal M}):=\{ \l\in \BR^{(m,n)}\,|\ {\Cal M}\cdot \l=0\ \}.$
The Eisenstein series $E_{k,{\Cal M}}^{(n)}$ is defined by
$$\align
E_{k,\M}^{(n)}(Z,W):=&\sum
\limits_{\left(\matrix A&B\\ C&D\endmatrix\right)\in
{\Cal R}}\,det\,(CZ+D)^{-k}\cdot e^{2\pi i\s(\M W(CZ+D)^{-1}C\,^tW)}\\
 &\cdot \sum\limits_{\lambda\in \Lambda}\, e^{2\pi i\s (\M((AZ+B)(CZ+D)^{-1}\,
^t\l+2\l\,^t(CZ+D)^{-1}\,^tW))},
\endalign$$
where $(Z,W)\in H_{n,m}.$ Now we assume that $k>n+m+1$ and $k$ is {\it even}.
Then according to [Z], Theorem 2.1, $E_{k,\M}^{(n)}(Z,W)$ is a nonvanishing
Jacobi form in $J_{k,\M}(\G_n).$ By Main Theorem, $(E_{k,\M}^{(n)})_{\tau}$
is a $\Hom(V_{\tau},\BC)$-valued modular form of type ${\hat {\tau}}\otimes
det^k$. We define the automorphic factor $j:Sp(n,\BR)\times H_n
\longrightarrow GL(n,\BC)$ by
$$j(g,Z):=cZ+d,\ \ \ \ g=
\left(\matrix a & b\\ c & d\endmatrix\right)\in Sp(n,\BR),\ \ \ Z\in H_n.$$
Then according to (5.9), for any homogeneous pluriharmonic polynomial $P$
with respect to $(2{\Cal M})^{-1},$ we obtain the following identity
$$\align
 & det\,j(M,Z)^k\,\sum\limits_{\gamma\in {\Cal R}}\sum\limits_{\lambda\in
\Lambda}\, det\,j(\gamma,Z)^{-k}\cdot
P(4\pi i{\Cal M}\l \,^tj(\gamma,Z)^{-1})
\cdot e^{2\pi i\s (\M\cdot \gamma<Z>\cdot ^t\l)}\\
=& \sum\limits_{\gamma\in {\Cal R}}\sum\limits_{\l\in \Lambda}\,
det\,j(\gamma,M<Z>)^{-k}\cdot P(4\pi i\M \l\,^tj(\gamma M,Z)^{-1})\cdot
e^{2\pi i\s (\M\cdot \gamma M<Z>\cdot ^t\l)}
\endalign$$
for all $M\in \Gamma_n$ and $Z\in H_n.$
\vskip 0.5cm
\noindent
{\smc Acknowledgements\,:}\ This work was done in part during my stay at
the Max-Planck-Institut f{\"u}r Mathematik. I am very grateful to the
institute for the hospitality and support.\par

\vskip 1cm {\smc Department\ of\ Mathematics\par Inha\
University\par  Incheon 402-751
\par
Republic of Korea} \vskip 0.3cm {\tt email\ address\,:\,\
jhyang\@inha.ac.kr}

\end